\documentclass[11pt]{amsart}
\usepackage{amssymb,amsfonts,amsthm}
\usepackage[all]{xy}

\newcommand{\IN}{\mathbb N}
\newcommand{\IR}{\mathbb R}
\newcommand{\IQ}{\mathbb Q}
\newcommand{\ID}{\mathbb D}
\newcommand{\II}{\mathbb I}

\newcommand{\I}{\mathcal I}
\newcommand{\wdi}{\mathrm{wd}}

\newcommand{\w}{\omega}
\newcommand{\C}{\mathcal C}
\newcommand{\hl}{\mathrm{hl}}

\newcommand{\U}{\mathcal{U}}

\newcommand{\F}{\mathcal{F}}
\newcommand{\Ra}{\Rightarrow}

\newcommand{\nw}{\mathrm{nw}}
\newcommand{\hd}{\mathrm{hd}}

\newcommand{\M}{\mathcal M}

\newcommand{\decc}{\overline{\mathrm{dec}}}

\newtheorem{thm}{Theorem}[section]
\newtheorem{theorem}[thm]{Theorem}
\newtheorem{lem}[thm]{Lemma}
\newtheorem{prop}[thm]{Proposition}
\newtheorem{proposition}[thm]{Proposition}
\newtheorem{cor}[thm]{Corollary}
\theoremstyle{definition}
\newtheorem{dfn}[thm]{Definition}
\newtheorem{ex}[thm]{Example}
\newtheorem{question}[thm]{Question}
\newtheorem{rem}[thm]{Remark}
\newtheorem{problem}[thm]{Problem}

\begin{document}
\title[On weak homeomorphisms]{Topological properties preserved by weakly discontinuous maps and weak homeomorphisms}
\author{Taras Banakh, Bogdan Bokalo, Nadiya Kolos}
\address{Jan Kochanowski University in Kielce (Poland) and Department of
Mathematics, Ivan Franko National University of Lviv (Ukraine)}
\email{t.o.banakh@gmail.com, b.m.bokalo@gmail.com, nadiya\_kolos@ukr.net}
\dedicatory{Dedicated to the 120th birthday of P.S.~Aleksandrov}
\begin{abstract} A map $f:X\to Y$ between topological
spaces is called {\em weakly discontinuous} if each
subspace $A\subset X$ contains an open dense subspace $U\subset A$ such that the restriction $f|U$ is continuous.
A bijective map $f:X\to Y$ between topological spaces is
called a {\em weak homeomorphism} if $f$ and $f^{-1}$ are
weakly discontinuous. We study properties of topological spaces
preserved by weakly discontinuous maps and weak homeomorphisms. In particular, we show that weak homeomorphisms preserve network weight, hereditary Lindel\"of number, dimension. Also we classify infinite zero-dimensional $\sigma$-Polish metrizable spaces up to a weak homeomorphism and prove that any such space $X$ is weakly homeomorphic to one of 9 spaces: $\w$, $2^\w$, $\IN^\w$, $\IQ$, $\IQ\oplus 2^\w$, $\IQ\times 2^\w$, $\IQ\oplus\IN^\w$, $(\IQ\times 2^\w)\oplus\IN^\w$,  $\IQ\times\IN^\w$.
\end{abstract}

\maketitle
\section{Introduction} \label{s1}

In this paper we detect topological properties preserved by weakly discontinuous maps and weak homeomorphisms.

By definition, a map $f:X\to Y$ between topological spaces is {\em
weakly discontinuous} if each subspace $A\subset X$ contains an open dense subspace $U\subset A$ such that the restriction $f|U$ is continuous.
Such maps were introduced by Vinokurov \cite{Vino} and studied in details in
\cite{AB, BB, BB1, BBK, BBK2, BK1, BK2, BK3, BM, KM, KS}.
Also they appear naturally in Analysis, see \cite{BKMM, CR, HOR}.


A bijective map $f:X\to Y$ between
topological spaces is called a {\em weak homeomorphism} if $f$
and $f^{-1}$ are weakly discontinuous. In this case we say that the topological spaces $X,Y$ are {\em weakly homeomorphic}. In particular, we show that
if $X,Y$ are weakly homeomorphic perfectly paracompact spaces,
then
\begin{enumerate}
\item $nw(X)=nw(Y)$;
\item $hd(X)=hd(Y)$;
\item $\dim X=\dim Y$;
\item $X$ is hereditarily Baire iff so is the space $Y$;
\item $X$ is analytic iff so is the space $Y$;
\item $X$ is $\sigma$-compact iff so is the space $Y$;
\item $X$ is $\sigma$-Polish iff so is the space $Y$.
\end{enumerate}
A topological space $X$ is called {\em $\sigma$-Polish} if it can be written as the countable union $X=\bigcup_{n\in\w}X_n$ of closed Polish subspaces.

In Sections~\ref{s2}--\ref{s4} we detect local and
global properties of topological spaces, preserved by weakly discontinuous maps and weak
homeomorphisms. In Section~\ref{s5} we classify zero-dimensional $\sigma$-Polish spaces up to weak homeomorphism and prove that each infinite zero-dimensional $\sigma$-Polish metrizable space $X$ is weakly homeomorphic to one of 9 spaces: $\w$, $2^\w$, $\IN^\w$, $\IQ$, $\IQ\oplus 2^\w$, $\IQ\times 2^\w$, $\IQ\oplus\IN^\w$, $(\IQ\times 2^\w)\oplus\IN^\w$, $\IQ\times\IN^\w$.

\subsection{Terminology and notations}
 Our terminology and notation are
standard and follow \cite{Ar} and \cite{En}. A ``space''
always means a ``topological space''. Maps between topological
spaces are not necessarily continuous.

By $\IR$ and $\IQ$ we denote the spaces of real and  rational
numbers, respectively; $\w$ stands for the space of finite ordinals
(= non-negative integers) endowed with the discrete topology. The set $\w\setminus\{0\}$ of finite positive ordinals (= natural numbers) is denoted by $\IN$. We
shall identify cardinals with the smallest ordinals of the given
cardinality.

For a subset $A$ of a space $X$ by 
$\bar{A}$ we denote the closure of $A$ in $X$.
For a function $f:X\to Y$
between topological spaces by $C(f)$ and $D(f)=X\setminus C(f)$ we
denote the sets of continuity and discontinuity points of $f$,
respectively.

Now we recall definitions of some cardinal invariants of topological spaces. For a
topological space $X$
\begin{itemize}
\item its {\em network weight} $\nw(X)$ is the smallest size $|\mathcal N|$
of a family $\mathcal N$ of subsets of $X$ such that for each point
$x\in X$ and each neighborhood $U\subset X$ of $x$ there is a set $N\in\mathcal N$ such that $x\in N\subset U$;
\item its {\em hereditary Lindel\"of number} $hl(X)$ is the smallest cardinal $\kappa$ such that each open cover of a subspace $Y\subset X$ has a subcover of cardinality $\le\kappa$;
\item its {\em hereditary density} $hd(X)$ is the smallest cardinal $\kappa$ such that each subspace $Y\subset X$ contains a dense subset of cardinality $\le\kappa$.
\end{itemize}

\section{Topological properties, preserved by weakly discontinuous maps}\label{s2}

In this section we discuss weakly discontinuous maps and detect topological properties preserved by such maps. We recall that a function $f:X\to Y$ between topological spaces is {\em weakly discontinuous} if any subset $A\subset X$ contains an open dense subset $U\subset A$ such that the restriction $f|U$ is continuous. Observe that a function $f:X\to Y$ is weakly discontinuous if and only if every non-empty subset $A\subset X$ contains a non-empty relatively open subset $V\subset A$ such that the restriction $f|V$ is continuous. This simple characterization implies the following useful (and known) fact.

\begin{lem}\label{l:comp} For two weakly discontinuous functions $f:X\to Y$ and $g:Y\to Z$ between topological spaces $X,Y,Z$ the composition $g\circ f:X\to Z$ is weakly discontinuous.
\end{lem}

\begin{proof} Given any non-empty subset $A\subset X$ use the weak discontinuity of $f$ and find a non-empty relatively open subset $U\subset A$ such that $f|U$ is continuous. Next, consider the non-empty set $B=f(U)$ and using the weak discontinuity of the function $g$, find a non-empty open set $V\subset B$ such that $f|V$ is continuous. The continuity of the function $f|U$ implies that the set $W=(f|U)^{-1}(V)$ is open in $U$. Then $W$ is a non-empty relatively open subset of $A$ and the restriction $g\circ f|W=(g|V)\circ (f|W)$ is continuous.
\end{proof}

Given a function $f:X\to Y$ between topological spaces and a subset $A\subset X$, let $D(f|A)$ be the set of discontinuity points of the restriction $f|A$ and $\bar D(f|A)$ be the closure of $D(f|A)$ in $X$. If $f$ is weakly discontinuous, then for any non-empty closed subset $A\subset X$ the sets $D(f|A)$ and $\bar D(f|A)$ are nowhere dense in $A$.

Put $\bar D_0(f)=X$ and for every ordinal $\alpha>0$ consider the closed subset $$\bar D_\alpha(f)=\bigcap_{\beta<\alpha}\bar D\big(f|\bar D_\beta(f)\big)$$ of $X$. The smallest ordinal $\alpha$ such that $\bar D_{\alpha+1}(f)=\bar D_\alpha(f)$ is called the {\em index of weak discontinuity} of $f$ and denoted by $\wdi(f)$. Since $\big(\bar D_\alpha(f)\big)_{\alpha<\wdi(f)}$ is a strictly decreasing transfinite sequence of closed sets in $X$, its length $\wdi(f)$ cannot exceed $\hl(X)^+$, the successor cardinal of the hereditary Lindel\"of number $\hl(X)$ of $X$ (i.e., $\wdi(f)<hl(X)^+$). It follows from the definition that a function $f:X\to Y$ is weakly discontinuous if and only if $\bar D_{\wdi(f)}(f)=\emptyset$.

In this case $X=\bigcup_{\alpha<\wdi(f)}\bar D_\alpha(f)\setminus \bar D_{\alpha+1}(f)$ and $f$ can be written as the composition $\tilde f\circ i$ of the identity function $$i:X\to\bigoplus_{\alpha<\wdi(f)}\bar D_{\alpha+1}(f)\setminus \bar D_\alpha(f), \quad i:x\mapsto x,$$ and the continuous function $\tilde f:\bigoplus_{\alpha<\wdi(f)}\bar D_{\alpha+1}(f)\setminus \bar D_\alpha(f)\to Y$, $\tilde f:x\mapsto f(x)$.

 Here by $\bigoplus_{\alpha<\wdi(f)}\bar D_{\alpha+1}(f)\setminus \bar D_\alpha(f)$ we denote the union $X=\bigcup_{\alpha<\wdi(f)}\bar D_{\alpha+1}(f)\setminus \bar D_\alpha(f)$ endowed with the topology of the topological sum of the family $\big(\bar D_{\alpha+1}(f)\setminus \bar D_\alpha(f)\big)_{\alpha<\wdi(f)}$. This topology consists of all subsets $U\subset X$ such that for every ordinal $\alpha<\wdi(f)$ the intersection $U\cap (\bar D_{\alpha+1}(f)\setminus \bar D_\alpha(f))$ is relatively open in $\bar D_{\alpha+1}(f)\setminus \bar D_\alpha(f)$. The following lemma implies that the function $i$ is weakly discontinuous.

\begin{lem}\label{l:sum} Let $\lambda$ be a non-zero ordinal and $(F_\alpha)_{\alpha\le \lambda}$ be a transfinite sequence of closed subsets of a topological space $X$ such that $F_0=X$, $F_\lambda=\emptyset$, $F_{\alpha+1}\subset F_\alpha$ for any $\alpha<\lambda$ and $F_\alpha=\bigcap_{\beta<\alpha}F_\beta$ for any limit ordinal $\alpha\le\lambda$. Then the identity map $$i:X\to\bigoplus_{\alpha<\lambda}F_\alpha\setminus F_{\alpha+1},\quad i:x\mapsto x,$$
is weakly discontinuous.
\end{lem}

\begin{proof} Given a non-empty subset $A\subset X$ we need to find a non-empty open set $U\subset A$ such that the restriction $i|U$ is continuous. Let $\beta$ be the smallest ordinal such that $A\subset F_\beta$. By the minimality of $\beta$, the relatively open subset $U=A\setminus F_{\beta+1}$ of $A$ is not empty. Then the restriction $$i|U:U\to F_\beta\setminus F_{\beta+1}\subset \bigoplus_{\alpha<\kappa}F_\alpha\setminus F_{\alpha+1}$$
is continuous.
\end{proof}

For a function $f:X\to Y$ between topological spaces, its {\em closed decomposition number}
$\decc(f)$ is defined as the smallest cardinality $|\mathcal{C}|$ of
a cover $\mathcal{C}$ of $X$ by closed or finite subsets such that $f|C$ is continuous
for each $C\in\mathcal{C}$. Observe that a function $f:X\to Y$ is
continuous iff $\decc(f)=1$.
Now we are going to give some upper bounds on the closed decomposition number $\decc(f)$ of a weakly discontinuous function $f$.

For a topological space $X$ the {\em large
pseudocharacter} $\Psi(X)$   is equal to the smallest cardinal
$\kappa$ such that each open set $U\subset X$ can be written as the union $U=\bigcup\F$ of a family $\mathcal{\F}$ consisting $\le\kappa$ many closed or finite subsets of $X$. It is easy to see that the large pseudocharacter $\Psi(X)$ of a regular space $X$ does not exceed its hereditary Lindel\"of number (i.e., $\Psi(X)\le \hl(X)$).



The following upper bound for the closed decomposition number of a weakly discontinuous map was given in Proposition 5.3 and Theorem 5.4 of \cite{BB}.

\begin{prop}\label{4.13} If\/ $f:X\to Y$ is a weakly discontinuous map between topological spaces, then $\decc(f)\le |\wdi(f)|\cdot \Psi(X)$. If the space $X$ is regular, then $|\wdi(f)|\cdot \Psi(X)\le \hl(X)$. If the space $X$ is paracompact, then $\decc(f)\le\Psi(X)$. If the space $X$ is perfectly paracompact, then $\decc(f)\le\w$.
\end{prop}

We recall that a topological space $X$ is {\em perfectly paracompact} if $X$ is paracompact and each open subset of $X$ is of type $F_\sigma$.

Now we detect some topological properties
preserved by weakly discontinuous maps. In the sequel we shall
identify a topological property with the class of topological spaces having that
property.

Let $\kappa$ be a cardinal. We shall say that a  property $\mathcal P$ of topological
spaces is
\begin{itemize}
\item {\em topological} if for any homeomorphic spaces $X,Y$ the space $X$ has property $\mathcal P$ if and only if $Y$ has that property;
\item {\em closed-hereditary} if for any space $X$ with property $\mathcal P$ every closed or finite subspace $Y$ of $X$ has property $\mathcal P$;
\item {\em open-hereditary} if for any space $X$ with property $\mathcal P$ every open subspace $Y$ of $X$ has property $\mathcal P$;
\item {\em projective} if for any continuous surjective map $f:X\to Y$ from a space $X$ with property $\mathcal P$ the image $Y=f(X)$ has property $\mathcal P$;
\item {\em $\kappa$-additive} if a space $X$ has property $\mathcal P$ whenever $X$ has a cover by $\le\kappa$ many closed subspaces with property $\mathcal P$;
\item {\em $\kappa$-summable} if for any family $\mathcal{C}$ of spaces with property $\mathcal P$ and $|\C|\le\kappa$ the topological sum $\oplus\C$ of $\C$ has property $\mathcal P$.
\end{itemize}


\begin{prop}\label{5.1} Let $f:X\to Y$ be a surjective weakly discontinuous map between topological spaces. If $X$ has some $\decc(f)$-summable closed-hereditary projective property $\mathcal P$, then the space $Y=f(X)$ has that property, too.
\end{prop}

\begin{proof} By the definition of the cardinal $\decc(f)$, we can find a
 cover $\mathcal{C}$ of $X$ by closed or finite subsets such that
$|\mathcal{C}|=\decc(f)$ and $f|C$ is continuous for
every $\mathcal{C}$. Consider the topological sum
$\oplus\mathcal{C}=\bigcup_{C\in\C}\{C\}\times C$ of the family $\mathcal{C}$. Taking into account that the property $\mathcal P$ is closed-hereditary and $|\mathcal C|$-summable, we
conclude that the topological sum $\oplus\mathcal{C}$ has property $\mathcal P$. Consider the
surjective continuous map $\tilde f:\oplus \mathcal{C}\to Y$ defined
by $\tilde f(C,x)=f(x)$ for any $C\in\mathcal{C}$ and $x\in C$. Then the space $Y$
belongs to $\mathcal P$, being a continuous image of the space
$\oplus\mathcal{C}$ possessing the projective property $\mathcal P$.
\end{proof}

\begin{prop}\label{5.1aa} Let $f:X\to Y$ be a surjective weakly discontinuous map between topological spaces. If $X$ has some $\wdi(f)$-summable closed-hereditary open-hereditary projective property $\mathcal P$, then the space $Y=f(X)$ has that property, too.
\end{prop}

\begin{proof} Consider the decreasing transfinite sequence $(\bar D_\alpha(f))_{\alpha<\wdi(f)}$ of closed sets in $X$ and observe that for every $\alpha<\wdi(f)$ the restriction $f|\bar D_\alpha(f)\setminus\bar D_{\alpha+1}(f)$ is continuous. Moreover, the closed and open heredity of the property $\mathcal P$ guarantees that the space $\bar D_\alpha(f)\setminus\bar D_{\alpha+1}(f)$ has property $\mathcal P$. By the $|\wdi(f)|$-summability of $\mathcal P$, the topological sum $$X_\oplus=\bigoplus_{\alpha<\wdi(f)}\bar D_\alpha(f)\setminus\bar D_{\alpha+1}(f)$$ has property $\mathcal P$. Then $Y$ has property $\mathcal P$, being a continuous image of the space $X_\oplus$ having the projective property $\mathcal P$.
\end{proof}

Applying Propositions~\ref{4.13} and \ref{5.1} to $\omega$-summable projective properties, we get

\begin{cor}\label{5.2} Let $f:X\to Y$ be a surjective weakly discontinuous map from a
perfectly paracompact space $X$ onto a topological space $Y$. If the space $X$ has an $\omega$-summable closed-hereditary projective property
$\mathcal P$, then $Y$ has that property, too.
\end{cor}

Applying Proposition~\ref{5.1aa} to the class of hereditarily Lindel\"of
$\sigma$-compact spaces, we get

\begin{cor}\label{5.3} Let $f:X\to Y$ be a weakly discontinuous surjective
map from a Hausdorff topological space $X$ to a topological space $Y$. If the space $X$ is hereditarily
Lindel\"of and $\sigma$-compact, then so is the space $Y$.
\end{cor}

Another corollary of Proposition~\ref{5.1aa} concerns analytic
spaces. A topological space $X$ is called {\em analytic} if $X$ is a
continuous image of a Polish (= separable completely metrizable)
space. Applying Proposition~\ref{5.1aa} to the class of
analytic spaces we get

\begin{cor}\label{5.4} A topological space $X$ is analytic if and
only if it is the image of a Polish space $P$ under a weakly discontinuous map $f:P\to X$.
\end{cor}

Now we detect some cardinal functions that respect weakly discontinuous maps. We define a  cardinal function $\varphi$ on a
 class $\mathcal T$ of topological spaces to be
\begin{itemize}
\item {\em topologically invariant} if $\varphi(X)=\varphi(Y)$ for any homeomorphic
spaces $X,Y\in\mathcal T$;
\item {\em closed-hereditary}  if every closed or finite subspace $Y$ of any space $X\in\mathcal T$ belongs to $\mathcal T$ and $\varphi(Y)\le\varphi(X)$;
\item {\em open-hereditary} if every open subspace $Y$ of any space $X\in\mathcal T$ belongs to $\mathcal T$ and $\varphi(Y)\le\varphi(X)$;
\item {\em projective} if $\varphi(f(X))\le\varphi(X)$ for
every  continuous map $f:X\to Y$ between spaces $X,Y\in\mathcal T$;
\item {\em $\wdi$-projective} if $\varphi(f(X))\le\varphi(X)$ for
every weakly discontinuous map $f:X\to Y$ between spaces $X,Y\in\mathcal T$;
\item {\em additive} if $\varphi(X)\le\sum_{C\in\mathcal \mathcal{C}}\varphi(C)$ for any space $X\in\mathcal T$ and a cover $\C\subset\mathcal T$ of $X$ by closed or finite subsets;
\item {\em summable} if $\varphi(\oplus\C)\le\sum_{C\in\mathcal \mathcal{C}}\varphi(C)$ for any family $\C\subset \mathcal T$ with $\oplus\C\in\mathcal T$;
\item {\em global} if $\varphi(D)\ge|D|$ for any discrete space $D\in\mathcal T$.
\end{itemize}

\begin{prop}\label{5.1a} Let $\varphi$ be a summable projective closed-hereditary (and open-hereditary)  cardinal function on the class of topological spaces. For any weakly discontinuous map $f:X\to Y$ between topological spaces, we get $\varphi(Y)\le\decc(f)\cdot \varphi(X)$ (and $\varphi(Y)\le\wdi(f)\cdot \varphi(X)\le\hl(X)\cdot \varphi(X)$~).
\end{prop}

\begin{proof}

Apply Proposition~\ref{5.1} (and \ref{5.1aa}) to the projective property $\mathcal P$ of a topological space $Z$ to have $\varphi(Z)\le \decc(f)\cdot \varphi(X)$ (and $\varphi(Z)\le \wdi(f)\cdot\varphi(X)$).
\end{proof}

\begin{thm}\label{n5.1} A global summable closed-hereditary  cardinal function $\varphi$ on
the class of regular spaces is $\wdi$-projective if and only if $\varphi$ is projective and $\varphi\ge \hl$.
\end{thm}

\begin{proof} Let $\varphi$ be a global summable closed-hereditary cardinal
function on the class of regular spaces. To prove the ``if'' part,
assume that $\varphi$ is projective and $\varphi\ge \hl$. Given a
surjective weakly discontinuous map $f:X\to Y$ between regular spaces, we can apply Propositions~\ref{5.1a} and \ref{4.13}  to conclude that
$$\varphi(Y)\le\decc(f)\cdot \varphi(X)\le \hl(X)\cdot \varphi(X)\le\varphi(X),$$
which means that $\varphi$ is $\wdi$-projective.
\smallskip

To prove the ``only if'' part, assume that $\varphi$ is $\wdi$-projective. Then it is projective (because each continuous map is
weakly discontinuous). It remains to prove that $\varphi\ge \hl$.
Assuming the converse, find a regular space $X$ with
$\varphi(X)<\hl(X)$. By \cite[2.9]{Juh},  $\hl(X)=\sup\{|Z|:Z\subset X$ is
scattered$\}$. So, we can find a scattered subspace
$Z\subset X$ of cardinality $|Z|>\varphi(X)$. Fix a bijective map $f:Z\to D$ of $Z$ onto
a discrete space $D$ and consider the topological sum $Y=(X\setminus
Z)\oplus D$. Consider the bijective map $\tilde f:X\to Y$ which is
identity on $X\setminus Z$ and coincides with $f$ on $Z$. We claim
that $f$ is weakly discontinuous. Given any non-empty subspace $A\subset X$
we should find a non-empty open set $U\subset A$ such that $\tilde f|U$ is continuous.
Consider the scattered subspace $C=Z\cap A$ of $A$.
If $C$ is not dense in $A$, then the relatively open set $U=A\setminus\bar C$ is a required open set such that the map  $\tilde f|U$ is continuous.
If $C$ is dense in $A$, then we can find an isolated point $x$ of the scattered space $C$ and choose an open subset $U\subset A$ of $A$ such that $U\cap C=\{x\}$. The density of $C$ in $A$ implies that $\{x\}=U\cap C$ is dense in $U$ and by the regularity, $\{x\}=U$. Then $U=\{x\}$ is open in $A$ and the restriction $\tilde f|U$ is trivially continuous.

Therefore the map $\tilde f:X\to Y$ is weakly discontinuous. The
$\wdi$-projectivity of $\varphi$ implies that
$\varphi(Y)\le\varphi(X)<|Z|$. Taking into account that
$Y=(X\setminus Z)\oplus D$ and $\varphi$ is closed-hereditary and
global, we arrive to a contradiction:
$|Z|=|D|\le\varphi(D)\le\varphi(Y)<|Z|.$
\end{proof}

Applying Proposition~\ref{5.1a} to some concrete hereditary projective cardinal functions we
get

\begin{cor}\label{5.5} If $f:X\to Y$ is a weakly discontinuous surjective map between
topological spaces, then
\begin{enumerate}
\item $nw(Y)\le nw(X)$;
\item $hl(Y)\le hl(X)$;
\item $hd(Y)\le\max\{hd(X),hl(X)\}$.
\end{enumerate}
\end{cor}

Because of an example of a regular space $X$ with $hd(X)<hl(X)$ (see
\cite{Tod}, \cite[3.12.7]{En} or \cite{Ro}), Theorem~\ref{n5.1}
implies that the hereditary density $hd$ is not
$\wdi$-projective, which means that there exists a weakly discontinuous map $f:X\to Y$ between regular spaces such that
$hd(Y)>hd(X)$. Under the Set Theoretic Assumption $\Diamond$ we can additionally assume that
$hd(X)=\aleph_0$.

\begin{ex} Under $\Diamond$ A.Ostaszewski \cite{Os} has constructed a regular space $X$ which is uncountable, compact, scattered, and hereditarily separable. Then any bijective map $f:X\to D$ to a discrete space $D$ is weakly discontinuous but $hd(D)=|D|=|X|>\aleph_0=hd(X)$. This means that the class of regular hereditarily separable spaces is not $\wdi$-projective under $\Diamond$.
\end{ex}

On the other hand, Todorcevic \cite{Tod} has constructed a model
of ZFC without $S$-spaces, that is, regular hereditarily separable
non-Lindel\"of spaces. In such models the class of regular
hereditarily separable spaces is $\wdi$-projective.

\begin{rem} Assume that no $S$-space exists. Then for each weakly discontinuous surjective map $f:X\to Y$ from a regular hereditarily separable $X$ the space $X$ is hereditarily Lindel\"of, and by Corollary~\ref{5.5} the image $Y=f(X)$ is hereditarily separable.
Consequently, the $\wdi$-projectivity of the class of regular hereditarily separable spaces is independent of the axioms of Set Theory.
\end{rem}

\section{Weak homeomorphisms}\label{s3}

In this section we introduce weak homeomorphisms and establish
their basic properties.

\begin{dfn}\label{6.1} Topological spaces $X$ and $Y$ are
called {\em weakly homeomorphic} if there is a bijective weakly discontinuous map $f:X\to Y$ with weakly discontinuous inverse $f^{-1}:Y\to X$. In this case the map $f$ is called a {\em weak homeomorphism}.
\end{dfn}

Corollaries~\ref{5.3}, \ref{5.4}, \ref{5.5} give the following list of properties preserved by weak homeomorphisms.

\begin{proposition}\label{p:wh-pres} If topological spaces $X,Y$ are weakly homeomorphic, then
\begin{enumerate}
\item $\nw(X)=\nw(Y)$;
\item $\hl(X)=\hl(Y)$;
\item $\hd(X)\cdot\hl(X)=\hd(Y)\cdot\hd(X)$;
\item $X$ is analytic if and only if $Y$ is analytic.
\end{enumerate}
\end{proposition}

For weak homeomorphisms we have the following result resembling the classical Cantor-Bernstein Theorem in Set Theory.

\begin{thm}\label{bernstein} Two topological spaces $X,Y$ are
weakly homeomorphic if each of them is homeomorphic to a closed subspace of the other space.
\end{thm}

\begin{proof} Assume that $X$ is homeomorphic to a closed subset $Y_1\subset Y$ and $Y$ is homeomorphic to a closed subset $X_1$ of $X$. Fix homeomorphisms $f:X\to Y_1$ and $g:Y\to X_1$. Let $X_0=X$, $Y_0=Y$ and inductively define subsets $X_{n+1}=g(Y_n)$ and $Y_{n+1}=f(X_n)$ for $n\ge0$. The inclusions $X_1\subset X$ and $Y_1\subset Y$ imply $Y_2=f(X_1)\subset f(X)=Y_1$ and $X_2=g(Y_1)\subset g(Y)=X_1$. Proceeding by induction, we can show that $X_{n+1}\subset X_n$ and $Y_{n+1}\subset Y_n$ for all $n\in\w$. Moreover, the sets $X_n$, $n\in\w$, are closed in $X$ and the sets $Y_n$, $n\in\w$, are closed in $Y$.

Consider the sets $X_\infty=\bigcap_{n\in\w}X_n$ and $Y_\infty=\bigcap_{n\in\w}Y_n$ and observe that $f(X_\infty)=f(\bigcap_{n\in\w} X_n)=\bigcap_{n\in\w}Y_{n+1}=Y_\infty$. Observe that the set $X_\infty$ is closed in $X$ and the set $Y_\infty$ is closed in $Y$.

Define a bijective function $h:X\to Y$ letting $h(x)=f(x)$ for $x\in X_\infty\cup\bigcup_{n\in\w}(X_{2n}\setminus X_{2n+1})$ and $h(x)=g^{-1}(x)$ for $x\in\bigcup_{n\in\w}X_{2n+1}\setminus X_{2n+2}$. By analogy with Lemma~\ref{l:sum}, it can be shown that bijective function $h:X\to Y$ is a weak homeomorphism.
\end{proof}

\begin{question}\label{6.6} Are two (regular) spaces weakly homeomorphic if each of them is weakly homeomorphic to a closed subspace of the other space?
\end{question}

\begin{rem}\label{6.7} Observe that the Baire space $\IN^\omega$ and the
Cantor cube $2^\omega$ embed into each other, but fail to be
weakly homeomorphic. The reason is that weak
homeomorphisms preserve the $\sigma$-compactness of hereditarily Lindel\"of spaces, see
Corollary~\ref{5.3}. This shows that the closedness is essential in
Theorem~\ref{bernstein}.
\end{rem}

Next, we show that weakly homeomorphic spaces can be decomposed into unions of closed
homeomorphic subspaces.

\begin{thm}\label{6.4} If $h:X\to Y$ is a bijective map between topological spaces,
then for some set $\I$ of cardinality $|\I|\le \decc(h)\cdot\decc(h^{-1})$ there is a cover  $\{X_i:i\in\I\}$ of $X$ by closed or finite subsets and a cover $\{Y_i:i\in\I\}$ of $Y$ by closed or finite subsets such that for every $i\in\I$ the restriction $h|X_i$ is a
homeomorphism of $X_i$ onto $Y_i$.
\end{thm}

\begin{proof}
By definition of the cardinal $\decc(h)$, we can find
a cover $\{X_\alpha:\alpha<\decc(h)\}$ of the space $X$ by closed or finite subsets such
that the restriction $h|X_\alpha$ is continuous for each
$\alpha<\decc(h)$. By analogy, for the function
$h^{-1}:Y\to X$ there exists a cover $\{Y_\beta:\beta<\decc(h^{-1})\}$ of
the space $Y$ by closed or finite subspaces such that the restriction $h^{-1}|Y_\beta$ is
continuous for each $\beta<\decc(h^{-1})$. Now let $\mathcal
I=\decc(h)\times\decc(h^{-1})$ and  for each $i=(\alpha,\beta)\in \mathcal I$
let $X_i=X_\alpha\cap h^{-1}(Y_\beta)$ and $Y_i=Y_\beta\cap
h(X_\alpha)$. It can be shown that the covers $\{X_i:i\in\I\}$ and $\{Y_i:i\in\I\}$
have the required property.
\end{proof}

Theorem~\ref{6.4} implies the following corollary.

\begin{cor}\label{5.1b} Let $f:X\to Y$ be a bijective map between topological spaces $X,Y$ and $\mathcal P$ be a $\kappa$-additive closed-hereditary topological property where $\kappa=\decc(f)\cdot\decc(f^{-1})$. The space $X$ has property $\mathcal P$ if and only if the space $Y$ has that property.
\end{cor}

Combining Corollary~\ref{5.1b} with Proposition~\ref{4.13}, we get the following two corollaries.

\begin{cor}\label{5.1wh} Let $f:X\to Y$ be a weak homeomorphism between topological spaces $X,Y$ and $\mathcal P$ be a $\kappa$-additive closed-hereditary topological property where $\kappa=\hl(X)\cdot\Psi(X)\cdot\hl(Y)\cdot\Psi(Y)$. The space $X$ has property $\mathcal P$ if and only if the space $Y$ has that property.
\end{cor}

\begin{cor}\label{5.1wh-pp} Let $f:X\to Y$ be a weak homeomorphism between perfectly paracompact  spaces $X,Y$. If the space $X$ has an $\w$-additive closed-hereditary property $\mathcal P$, then the space $Y$ has that property, too.
\end{cor}

Applying Corollary~\ref{5.1wh} to some concrete topological properties, we
get

\begin{cor}\label{10.10} If $X,Y$ are weakly homeomorphic perfectly paracompact  spaces, then
\begin{enumerate}
\item $\nw(X)=\nw(Y)$;
\item $\hl(X)=\hl(Y)$;
\item $\hd(X)=\hd(X)$;
\item $\dim X=\dim Y$;
\item $X$ is analytic iff so is the space $Y$;
\item $X$ is $\sigma$-compact iff so is the space $Y$;
\item $X$ is $\sigma$-Polish iff so is the space $Y$.
\end{enumerate}
\end{cor}

A cardinal function $\varphi$ on the class of (regular) spaces is
defined to be {\em invariant under weak homeomorphisms} if $\varphi(X)=\varphi(Y)$
for any two weakly homeomorphic (regular) spaces $X,Y$.

The following theorem characterizing such cardinal functions is a counterpart of Theorem~\ref{n5.1}.

\begin{thm}\label{9.1} A global additive closed-hereditary cardinal function $\varphi$ on the class of regular spaces is invariant under weak homeomorphisms if and only if $\varphi$ is topologically invariant and $\varphi\ge \hl$.
\end{thm}

\begin{proof} To prove the ''if'' part, assume that $\varphi$ is topologically invariant  and $\varphi\ge \hl$. It suffices to show that $\varphi(Y)\le\varphi(X)$ for any two weakly homeomorphic regular spaces $X,Y$. By Corollary~\ref{5.5}, $\hl(Y)\le\hl(X)$. Let $\kappa=\varphi(X)$ and $\mathcal P$ be the class of regular spaces $Z$ with $\varphi(Z)\le\kappa$. The closed-heredity and additivity of $\varphi$ implies the closed-heredity and $\kappa$-additivity of the property $\mathcal P$. Since $\Psi(Y)\le\hl(Y)\le \hl(X)=\hl(X)\cdot \Psi(X)\le\varphi(X)=\kappa$, we may apply Corollary~\ref{5.1wh} to conclude that $Y$ has property $\mathcal P$ and hence $\varphi(Y)\le \kappa=\varphi(X)$.
\smallskip

To prove the ``only if'' part, assume that a global
closed-hereditary cardinal function $\varphi$ is
invariant under weak homeomorphisms. It is clear that $\varphi$ is
topologically invariant. So it remains to show that $\varphi(X)\ge
\hl(X)$ for each regular space $X$. Assuming the converse and
repeating the argument of the proof of Theorem~\ref{n5.1}, we can
construct a weak homeomorphism $h:X\to Y$ of $X$ onto a regular
space $Y$ containing a closed discrete subspace $D$ of size
$|D|>\varphi(X)$, which leads to a contradiction:
 $\varphi(X)<|D|\le\varphi(D)\le\varphi(Y)=\varphi(X)$.
\end{proof}

\begin{rem} As an example of an exotic cardinal function, invariant under weak homeomorphisms between regular spaces, let us consider the cardinal function $ccw(X)$ called the {\em closed covering weight} of $X$ and equal to the smallest cardinal $\kappa$ for which there is a cover of $X$ by $\le\kappa$ closed subspaces with weight $\le\kappa$. It is easy to see that $nw(X)\le ccw(X)\le w(X)$ for every topological space $X$. Both the inequalities can be strict: $ccw(X)<w(X)$ for any countable space $X$ with uncountable weight.

To construct a space $X$ with $nw(X)<ccw(X)$, take any space $Z$
with $\aleph_0=nw(Z)<w(Z)$ containing a dense countable discrete
subspace $D$. The countable power $X=Z^\w$ is a Baire space with
countable network weight. For every countable closed cover
$\mathcal{C}$ of $X$ the Baire Theorem gives a set $C\in\mathcal{C}$
with non-empty interior in $X$. Then $w(C)=w(Z)>\aleph_0$, which
means that $ccw(X)>\aleph_0$.
\end{rem}

By Corollary~\ref{10.10}, weak homeomorphisms between perfectly paracompact spaces preserve $\sigma$-compact spaces and $\sigma$-Polish spaces.
We shall show that they also preserve Polish spaces. To prove this fact, we need the following  continuity property of weak homeomorphisms.

\begin{thm}\label{loc0} Let $h:X\to Y$ be a weak homeomorphism between topological spaces. Then each non-empty closed subset $A\subset X$ contains a non-empty open subset $U\subset A$ such that $h(U)$ is open in its closure $\overline{h(U)}$ in $Y$ and $h|U:U\to h(U)$ is a homeomorphism.
\end{thm}

\begin{proof} Since $h$ is weakly discontinuous, the open subset $V=A\setminus\bar D(f|A)$ is dense in $A$. Let $B=\overline{h(V)}$ be the closure of $h(V)$ in $Y$. By the weak discontinuity of the map $h^{-1}$, there exists an open dense subset
$W\subset B$ such that $h^{-1}|W$ is continuous.  The continuity of
$h|V$ and density of $h(V)$ in $B$ implies that the set
$U=(h|V)^{-1}(W)=V\cap h^{-1}(W)$ is non-empty and open in $V$ and
hence in $A$.

Observe that $h(U)=h(V\cap h^{-1}(W))=h(V)\cap W=(h^{-1}|W)^{-1}(V)$
is open in $W$ as the preimage of the open set $V$ under the
continuous map $h^{-1}|W$. Since $W$ is open in $B$, the set $h(U)$
is open in $B$ and hence in $\overline{h(U)}\subset B$. Also
$h|U:U\to h(U)$ is a homeomorphism because the inverse map
$h^{-1}|h(U)$ is continuous being the restriction of the continuous
map $h^{-1}|W$ to the set $h(U)=h(V)\cap W$.
\end{proof}

Using Theorem~\ref{loc0} we shall prove that weak homeomorphisms preserve hereditarily Baire spaces. We recall that a topological space $X$ is {\em hereditarily Baire} if each closed subspace of $X$ is Baire. A space $X$ is {\em Baire} if the intersection $\bigcap\U$ of any countable family of open dense subsets is dense in $X$.

\begin{thm}\label{t:herBaire} Let $h:X\to Y$ be a weak homeomorphism between topological spaces. The space $X$ is hereditarily Baire if and only if $Y$ is hereditarily Baire.
\end{thm}

\begin{proof} To prove the ``if'' part, assume that $Y$ is hereditarily Baire. To derive a contradiction, assume that $X$ is not hereditarily Baire. Then we can find a non-empty closed subspace $Z\subset X$ which is not Baire. Replacing $Z$ by a suitable closed subspace, we can assume that  every non-empty open subspace of $Z$ is not Baire. By Theorem~\ref{4.13}, the space $Z$ contains a non-empty open subspace $U$ such that the restriction $h|U:U\to h(U)$ is a homeomorphism and $h(U)$ is open in its closure $F=\overline{h(U)}$ in $Y$. Since the space $Y$ is hereditarily Baire, the closed subspace $F$ of $Y$ is Baire and so is its open dense subspace $h(U)$ and the topological copy $U$ of $h(U)$. But this contradicts the choice of $Z$.
\end{proof}

\begin{rem} In general, weak homeomorphisms do not preserve Baire spaces. It is easy to construct two metrizable separable spaces $X=A\cup B$ and $X'=A'\cup B'$ such that $A,A'$ are homeomorphic to the space $\IQ$ of rational numbers, $B,B'$ are discrete spaces, $A$ is nowhere dense in $X$ and $A'$ is open and dense in $X$. It is easy to see that the spaces $X$ and $X'$ are weakly homeomorphic, $X$ is Baire but $X'$ is meager.
\end{rem}

We recall that a topological space $X$ is {\em Polish} if it is separable and its topology is generated by a complete metric. The Baire Theorem implies that each Polish space is hereditarily Baire. The converse is true for coanalytic spaces by the classical Hurewicz Theorem 21.18 \cite{Ke}.

\begin{cor}\label{polish} A metrizable space $X$ is Polish if and only if it is weakly homeomorphic to a Polish space.
\end{cor}

\begin{proof} The ``only if'' part is trivial. To prove the ``if'' part, assume that $X$ is weakly homeomorphic to a Polish space $Y$. By Corollary~\ref{10.10} the space $X$ is $\sigma$-Polish and by Theorem~\ref{t:herBaire}, $X$ is hereditarily Baire. By Hurewicz Theorem~21.18 \cite{Ke}, each coanalytic hereditarily Baire space is Polish, which implies that $X$ is Polish.
\end{proof}

The metrizability is not preserved by weak homeomorphisms
(because each scattered space is weakly homeomorphic to a
discrete space). However we have a partial result for perfectly normal compact
spaces. We recall that a normal space $X$ is {\em perfectly normal} if each open subset of $X$ is of type $F_\sigma$. It is well-known (and easy to see) that a compact space is perfectly normal if and only if it is hereditarily Lindel\"of.

\begin{prop}\label{10.15} A perfectly normal compact space $X$ is metrizable
if and only if $X$ is weakly homeomorphic to a metrizable space.
\end{prop}

\begin{proof} Assume that a perfectly normal compact space $X$ admits
a weak homeomorphism $h:X\to Y$ onto a metrizable space $Y$.
Applying Proposition~\ref{p:wh-pres}, we get $\hl(Y)\le \hl(X)\le\aleph_0$.
Then the metrizable space $Y$, being hereditarily Lindel\"of, is
separable and thus has countable network weight. Applying
Proposition~\ref{p:wh-pres}, we conclude that the compact space $X$ has
countable network weight and hence is metrizable according to
\cite[3.1.19]{En}.
\end{proof}

The perfect normality in this theorem is essential since each
scattered compact space is weakly homeomorphic to a discrete
space. On the other hand, a metrizable space that is weakly
homeomorphic to a  compact spaces, need not have countable base:
each discrete space is metrizable and weakly homeomorphic to
its Aleksandrov compactification.

Another property preserved by scattered homeomorphisms is the
$\mathcal{C}$-universality. We define a space $X$ to be (everywhere)
{\em $\mathcal{C}$-universal} for a class $\mathcal{C}$ of spaces if
for every space
$C\in\mathcal{C}$ (and every non-empty open subset $U\subset X$) there exists a closed embedding $e:C\to X$ (with
$e(C)\subset U$). For example, the Hilbert cube $[0,1]^\w$ is everywhere
$\M_0$-universal for the class $\M_0$ of compact metrizable spaces
while the Hilbert space $l^2$ is everywhere $\M_1$-universal for the
class $\M_1$ of Polish spaces (see, e.g. \cite{BRZ}).

\begin{prop}\label{10.11} Assume that a class $\mathcal{C}$ of topological spaces contains a Baire perfectly paracompact everywhere $\mathcal{C}$-universal space $U$. A perfectly paracompact space $X$ is $\mathcal{C}$-universal if and only if $X$ is weakly homeomorphic to a $\mathcal{C}$-universal space $Y$.
\end{prop}

\begin{proof} Let $h:X\to Y$ be a weak homeomorphism. Since the space $Y$ is $\mathcal{C}$-universal, it contains a closed subspace $Z$ homeomorphic to $U$. Using the perfect paracompactness of $X$ and Proposition~\ref{4.13}, find a countable closed cover $\{X_n:n\in\w\}$ of $X$ such that the restrictions $h|X_n$ are continuous for all $n\in\w$. By the same reason, the space $Z$ has a countable closed cover $\{Z_m:m\in\w\}$ such that for every $m\in\w$ the restriction $h^{-1}|Z_m$ is continuous. Then for every $n\in\w$ the set $Z_m^n=(h^{-1}|Z_m)^{-1}(X_n)=Z_m\cap h(X_n)$ is closed in $Z_m$ and is homeomorphic to the closed subset $h^{-1}(Z_m^n)=(h|X_n)^{-1}(Z_m)$ of $X_n$. Since the space $Z=\bigcup_{n,m\in\w}Z_m^n$ is Baire,  for some $m,n\in\w$ the set $Z_{m,n}$ has non-empty interior in $Z$ and hence $Z_m^n$ is $\mathcal{C}$-universal (because $Z$ is everywhere $\mathcal{C}$-universal). Then $X$ is $\mathcal{C}$-universal as well because it contains a closed topological copy $h(Z^n_m)$ of the $\mathcal{C}$-universal space $Z_m^n$.
\end{proof}

\section{Detecting local topological properties preserved by weak homeomorphisms}
\label{s4}

In this section we characterize local topological properties, preserved by weak homeomorphisms.

A property $\mathcal P$ of a topological space is said to
be
\begin{itemize}
\item {\em local} if a topological space $X$ has the property $\mathcal P$ if and only if each point $x\in X$ has a neighborhood $U$ with property $\mathcal P$;
\item {\em closed+open additive} if a topological space $X$ has the property $\mathcal P$ provided $X$ contains an open subset $U$ such that $U$ and $X\setminus U$ have the property $\mathcal P$;
\item {\em scattered}  if a topological space $X$ has property $\mathcal P$ if and only if each non-empty closed subspace $A\subset X$ contains a non-empty relatively open subset $U$ with property $\mathcal P$;
\item {\em preserved by weak homeomorphisms} if for any weak homeomorphism $h:X\to Y$ between topological spaces the space $X$ has property $\mathcal P$ if and only if $Y$ has $\mathcal P$.
\end{itemize}

The main result of this section is the following characterization theorem.

\begin{thm}\label{loc2} For a closed-hereditary open-hereditary topological property $\mathcal P$ the following conditions are equivalent:
\begin{enumerate}
\item $\mathcal P$ is scattered;
\item $\mathcal P$ is local and open+closed additive;
\item $\mathcal P$ is local and is preserved by weak homeomorphisms.
\end{enumerate}
\end{thm}

\begin{proof} We identify $\mathcal P$ with the class of topological spaces that possess the property $\mathcal P$.
\smallskip

$(1)\Ra(2)$ Assume that the property $\mathcal P$ is scattered. To prove that $\mathcal P$ is local, assume that $X$ has a cover $\U$ by open subsets with the property $\mathcal P$. To show that $X\in\mathcal P$ it suffices to find in each non-empty closed subset  $F\subset X$ a non-empty open subset $V\subset F$ with property $\mathcal P$. Find an open set $U\in\U$ such that $F\cap U\ne\emptyset$ and observe that the non-empty set $V=F\cap U$ is open in $F$, closed in $U$ and hence has the property $\mathcal P$ as $\mathcal P$ is closed-hereditary.

To show that $\mathcal P$ is open+closed additive, assume that a topological space $X$ contains an open set $U$ such that $U$ and $X\setminus U$ belong to $\mathcal P$. Since the property $\mathcal P$ is scattered, the inclusion $X\in\mathcal P$ will follow as soon as we prove that each non-empty closed subspace $F\subset X$ contains a non-empty open subspace $V$ with property $\mathcal P$. If $F\subset X\setminus U$, then $F$ has property $\mathcal P$ since this property is closed-hereditary. So, we assume that $F\not\subset X\setminus U$ and hence the open subset $V=F\cap U$ is non-empty and has the property $\mathcal P$ since $U\in\mathcal P$ and $V=F\cap U$ is closed in $U$.
\smallskip

$(2)\Ra(3)$ Now assume that the property $\mathcal P$ is local and open+closed additive. Let $f:X\to Y$ be a weak homeomorphism between topological spaces. Assume that $Y\in\mathcal P$ but $X\notin\mathcal P$. Let $\U$ be the family of all open subspaces of $X$ with property $\mathcal P$. The locality of $\mathcal P$ implies that the open set $U=\bigcup\U$ has property $\mathcal P$. Since the property $\mathcal P$ is open+closed additive, the closed set $F=X\setminus U$ does not have $\mathcal P$. In particular, $F$ is not empty. We claim that each non-empty open set $V\subset F$ does not have the property $\mathcal P$. Assuming the converse, and taking into account that $\mathcal P$ is open+closed additive, we conclude that the open subspace  $U\cup V$ of $X$ has property $\mathcal P$ and hence $U\cup V\in\U$, which is not possible as $U\cup V\not\subset\bigcup\U$. So, the closed set $F$ nowhere has property $\mathcal P$. By Theorem~\ref{loc0}, the space $F$ contains a non-empty open set $U$ such that $f(U)$ is open in its closure $\overline{f(U)}$ in $Y$ and $f|U:U\to f(U)$ is a homeomorphism. Since the property $\mathcal P$ is closed-hereditary, the closed subset $\overline{f(U)}$ of the space $Y\in\mathcal P$ has property $\mathcal P$. Since $\mathcal P$ is open-hereditary, the open subset $f(U)$ of $\overline{f(U)}$ has property $\mathcal P$. Since $\mathcal P$ is topological, the open set $U\subset F$ has property $\mathcal P$, which is a desired contradiction showing that $X\in\mathcal P$.
\smallskip

$(3)\Ra(1)$ Assume that the closed-hereditary property $\mathcal P$ is local and is preserved by weak  homeomorphisms.  To
show that $\mathcal P$ is scattered, we should check that a topological
space $X$ has property $\mathcal P$ provided each non-empty closed
subspace $A\subset X$ contains a non-empty open subspace $U\subset
A$  with $\mathcal P$. Using the latter property of $X$ and the
locality of $\mathcal P$, we may construct a transfinite
sequence $(X_\alpha)_{\alpha\le \lambda}$ of closed subsets of $X$
such that $X_0=X$, $X_\lambda=\emptyset$, $X_\beta=\bigcap_{\alpha<\beta}X_\alpha$ for any limit ordinal $\beta\le\lambda$ and for every ordinal $\alpha<\lambda$ the set
$X_\alpha\setminus X_{\alpha+1}$  is dense in $X_\alpha$ and has property $\mathcal P$. The
locality of $\mathcal P$ implies that the topological sum
$Y=\bigoplus_{\alpha<\beta}X_\alpha\setminus X_{\alpha+1}$ also has
property $\mathcal P$. Observe that the ``identity'' map $i:Y\to X$
is continuous while its inverse $i^{-1}:X\to Y$ is weakly discontinuous. Since $\mathcal P$ is preserved by weak homeomorphisms, the space $X$ belongs to $\mathcal P$, being weakly
homeomorphic to the space $Y\in\mathcal P$.
\end{proof}

Each property $\mathcal P$ of topological spaces
induces a scattered property called the $\mathcal P$-scatteredness.
Namely, we say that a regular space $X$ is {\em $\mathcal
P$-scattered} if each closed non-empty subspace $A\subset X$
contains a non-empty relatively open subspace $U\subset A$ with
property $\mathcal P$, see \cite{CD}.

Let us note the scatteredness is just $\mathcal P$-scatteredness for
the property $\mathcal P$ of being a singleton.

The $\mathcal P$-scatteredness can be
characterized as follows.

\begin{thm}\label{loc3} Let $\mathcal P$ be a local closed-hereditary open-hereditary topological property. For a (regular) topological space $X$ the following conditions are
equivalent:
\begin{enumerate}
\item $X$ is $\mathcal P$-scattered;
\item $X$ is weakly homeomorphic to a (regular) $\mathcal P$-scattered space $Y$;
\item there is a bijective continuous map $h:Y\to X$ from a (regular) space $Y$ possessing the property $\mathcal P$ whose inverse $h^{-1}$ is weakly discontinuous.
\end{enumerate}
\end{thm}

\begin{proof} It suffices to prove the implications $(1)\Ra(3)\Ra(2)\Ra(1)$.

The proof $(1)\Ra(3)$ repeats the proof of the implication $(3)\Ra(1)$ of
Theorem~\ref{loc2}, $(3)\Ra(2)$ is trivial, and $(2)\Ra(1)$ follows
from Theorem~\ref{loc2} (applied to the property of being a
$\mathcal P$-scattered space).
\end{proof}

If $\mathcal P$ is the class of locally compact spaces, then the
$\mathcal P$-scatteredness is referred to as $k$-scatteredness. More
precisely, a regular space $X$ is called {\em $k$-scattered} if each
non-empty closed subspace $A\subset X$ contains a compact subspace
$K\subset A$ with non-empty interior in $A$, see \cite{CD}.

\begin{cor}\label{7.1}  For a (regular) topological space $X$ the following conditions are
equivalent:
\begin{enumerate}
\item $X$ is $k$-scattered;
\item $X$ is weakly homeomorphic to a (regular) $k$-scattered space $Y$;
\item there is a bijective continuous map $h:Y\to X$ from a (regular) locally compact space $Y$ whose inverse $h^{-1}$ is weakly discontinuous;
\item $X$ is weakly homeomorphic to a compact (Hausdorff) space;
\end{enumerate}
\end{cor}

\begin{proof} The equivalence of the first three conditions follows from Theorem~\ref{loc3} while $(4)\Ra(2)$ is trivial. The proof will be complete if we prove the implication $(3)\Ra(4)$.

Assume that $h:X\to Y$ is a weak homeomorphism of $X$ onto a
locally compact (regular) space $Y$. Fix any point $y_\infty\in Y$ and
consider a new (regular) topology $\tau^*$ on $Y$ coinciding with the original
topology at each point $y\ne y_\infty$ and such that a neighborhood
base at the point $y_\infty$ consists of the sets $O(y_\infty)\cup
(Y\setminus K)$, where $O(y_\infty)$ is a neighborhood of $y_\infty$
in $Y$ and $K$ is a compact set in $Y$. It is easy
to see that $(Y,\tau^*)$ is a compact (Hausdorff) space and $h:X\to (Y,\tau^*)$
is a weak homeomorphism.
\end{proof}

Finally we shall prove a useful decomposition lemma, which will be used in the proof of the classification Theorem~\ref{t:class}.

\begin{lem}\label{l:decomp} Let $\mathcal P_n$, $n\in\IN$, be local topological properties. Every (metrizable) topological space $X$ is weakly homeomorphic to the topological sum $\bigoplus_{n\in\w}X_n$ of spaces such that
\begin{itemize}
\item for every $n\in\IN$ the space $X_n$ has property $\mathcal P_n$ (and is metrizable);
\item the space $X_0$ is homeomorphic to a closed subset of $X$ and every non-empty subset $U\subset X$ has properties $\mathcal P_n$ for no $n\in\w$.
\end{itemize}
\end{lem}

\begin{proof} It is well-known that every ordinal $\alpha$ can be uniquely written as $\alpha=\beta+n$ where $\beta$ is a limit ordinal and $n\in\w$ is finite. The number $n$ will be called the {\em integer part} of $\alpha$ and will be denoted by $\lfloor\alpha\rfloor$. An ordinal $\alpha$ is called {\em odd} (resp. {\em even}) if its integer part $\lfloor\alpha\rfloor$ is odd (resp. even). Observe that each odd ordinal is not limit. For an ordinal $\alpha$ with positive integer part by $\alpha-1$ we denote the unique ordinal such that $\alpha=(\alpha-1)+1$.

Let $X$ be a topological space. For every $n\in\IN$ and every subspace $A\subset X$ let $U_n(A)$ be the union of all open subsets $U\subset A$ that have the property $\mathcal P_n$. Since the class $\mathcal P_n$ is local, the space $U_n(A)$ belongs to $\mathcal P_n$.

Let $Z_0=X$ and for every ordinal $\alpha$ define a closed subset $Z_\alpha$ of $X$ by the recursive formula
$$Z_\alpha=\begin{cases}\bigcap_{\beta<\alpha}Z_\beta&\mbox{if $\alpha$ is limit}\\
Z_{\alpha-1}\setminus U_{\lfloor \alpha\rfloor}(Z_{\alpha-1})&\mbox{otherwise}.
\end{cases}
$$
Since $(Z_\alpha)_\alpha$ is a decreasing sequence of closed subsets of $X$, there is a limit ordinal $\lambda$ such that $Z_{\lambda+n}=Z_\lambda$ for all $n\in\IN$, which means that  $U_n(Z_\lambda)=\emptyset$ for all $n\in\w$ and hence no non-empty open set of $Z_\lambda$ has a property $\mathcal P_n$ for some $n\in\w$. It is clear that $X=Z_\lambda\cup\bigcup_{\alpha<\lambda}Z_{\alpha}\setminus Z_{\alpha+1}$ and $X$ is weakly homeomorphic to the topological sum $Z_\lambda\oplus\bigoplus_{\alpha<\lambda}Z_{\alpha}\setminus Z_{\alpha+1}$.

Put $X_0=Z_\lambda$ and for every $n\in\IN$ let $$X_n=\textstyle{\bigoplus}\{Z_{\alpha-1}\setminus Z_\alpha:\alpha<\lambda,\;\lfloor\alpha\rfloor=n\}=
\textstyle{\bigoplus}\{U_n(Z_{\alpha-1}):\alpha<\lambda,\;\lfloor\alpha\rfloor=n\}.
$$
The locality of the class $\mathcal P_n$ guarantees that $X_n\in\mathcal P_n$.
It is clear that the topological sum $\bigoplus_{n\in\w}X_n$ is equal to the topological sum $Z_\lambda\oplus\bigoplus_{\alpha<\kappa}Z_{\alpha}\setminus Z_{\alpha+1}$. So, $X$ is weakly homeomorphic to $\bigoplus_{n\in\w}X_n$.
\end{proof}

\section{Classifying zero-dimensional $\sigma$-Polish spaces up to a weak homomorphism}\label{s5}

In this section we shall classify infinite zero-dimensional metrizable $\sigma$-Polish spaces up to a weak homeomorphisms. We recall that  a topological space $X$ is {\em $\sigma$-Polish} if $X$ can be written as the countable union $X=\bigcup_{n\in\w}X_n$ of closed Polish subspaces of $X$. So, any Polish or metrizable $\sigma$-compact space is $\sigma$-Polish. In particular, the space $\IQ$ of rational numbers is $\sigma$-Polish and so are the countable powers $2^\w$ and $\IN^\w$ of the doubleton $2=\{0,1\}$ and the discrete space $\IN$ of natural numbers. It is clear that the class of $\sigma$-Polish spaces is closed under finite products and countable topological sums.
A topological space $X$ is called {\em zero-dimensional} if $X$ has a base of topology consisting of open-and-closed subsets.

In Theorem~\ref{t:class} we shall prove that every infinite zero-dimensional $\sigma$-Polish metrizable space is homeomorphic to one of 9 spaces included in the following diagram (in which an arrow $X\to Y$ between two spaces $X,Y$ indicates that $X$ is homeomorphic to a closed subspace of $Y$):
$$
\xymatrix{
&2^\w\ar[r]\ar[rd]&\IN^\w\ar[r]&\IQ\oplus\IN^\w\ar[rd]\ar[r]&\IN^\w\\
\w\ar[r]&\IQ\ar[r]&\IQ\oplus 2^\w\ar[r]\ar[ru]&\IQ\times 2^\w\ar[r]&(\IQ\times 2^\w)\oplus \IN^\w\ar[u]
}
$$

We shall need the following known topological characterizations of the spaces $\IQ$, $2^\w$,  $\IN^\w$, $\IQ\times 2^\w$ and $\IQ\times\IN^\w$ due to Sierpi\'nski \cite{S}, Brouwer \cite{Brouwer}, Alexandroff, Urysohn \cite{AU}, and van Mill \cite{vM81}.

\begin{prop}\label{p:char} A metrizable zero-dimensional space $X$ is homeomorphic to
\begin{enumerate}
\item $\IQ$ if and only if $X$ is countable and has no isolated points;
\item $2^\w$ if and only if $X$ is compact and has no isolated points;
\item $\IN^\w$ if and only if $X$ is Polish and nowhere locally compact;
\item $\IQ\times 2^\w$ if and only if $X$ is $\sigma$-compact, nowhere countable and nowhere locally compact;
\item $\IQ\times\IN^\w$ if and only if $X$ is $\sigma$-Polish, nowhere $\sigma$-compact and nowhere Polish.
\end{enumerate}
\end{prop}

These characterizations imply the following known embedding results.

\begin{prop}\label{p:emb} A metrizable zero-dimensional space $X$ is homeomorphic to a closed subspace of
\begin{enumerate}
\item $\IQ$ if and only if $X$ is countable;
\item $2^\w$ if and only if $X$ is compact;
\item $\IN^\w$ if and only if $X$ is Polish;
\item $\IQ\times 2^\w$ if and only if $X$ is $\sigma$-compact;
\item $\IQ\times\IN^\w$ if and only if $X$ is $\sigma$-Polish.
\end{enumerate}
\end{prop}

We shall also need the following three tests due to Aleksandrov \cite{A16} and Hurewicz \cite{Hur} (see also \cite[29.1, 21.18, 21.19]{Ke}).

\begin{prop}\label{p:hur} A Borel subset $X$ of a Polish space is
\begin{enumerate}
\item uncountable if and only if $X$ contains a subspace homeomorphic to $2^\w$;
\item not Polish if and only if $X$ contains a closed subspace homeomorphic to $\IQ$;
\item not  $\sigma$-compact if and only if $X$ contains a closed subspace homeomorphic to $\IN^\w$.
\end{enumerate}
\end{prop}

The following classification theorem is the main result of this section.

\begin{thm}\label{t:class} Let $X$ be an infinite zero-dimensional metrizable space.
\begin{enumerate}
\item If $X$ is Polish, then $X$ is weakly homeomorphic to one of 3 spaces: $\w$, $2^\w$, $\IN^\w$.
\item If $X$ is $\sigma$-compact, then $X$ is weakly homeomorphic to one of 5 spaces: $\w$, $2^\w$, $\IQ$, $\IQ\oplus 2^\w$, $\IQ\times 2^\w$.
\item If $X$ is $\sigma$-Polish, then $X$ is weakly homeomorphic to one of 9 spaces: $\w$, $2^\w$, $\IN^\w$, $\IQ$, $\IQ\oplus 2^\w$, $\IQ\oplus\IN^\w$, $\IQ\times 2^\w$, $(\IQ\times 2^\w)\oplus\IN^\w$, $\IQ\times\IN^\w$.
\end{enumerate}
\end{thm}

\begin{proof} 1. First we assume that the space $X$ is Polish. This case has three subcases.
\smallskip

1.1. The space $X$ is countable. Then each closed subspace of $X$, being Polish and countable, has an isolated point, which implies that $X$ is scattered and hence is weakly homeomorphic to the discrete space $|X|=\w$.
\smallskip

1.2. The space $X$ is uncountable and $\sigma$-compact. In this case we shall prove that $X$ is $k$-scattered. Write $X$ as the countable union $X=\bigcup_{n\in\w}K_n$ of compact subsets. Given a non-empty closed subset $A\subset X$ we can apply the Baire Theorem and find $n\in\w$ such that the compact set $K=A\cap K_n$ has non-empty interior in $A$, witnessing that $X$ is $k$-scattered. By Corollary~\ref{7.1}, the $k$-scattered space $X$ is weakly homeomorphic to an uncountable compact Hausdorff space $K$.   By Corollary~\ref{5.5}, the compact space $K$ has network weight $\nw(K)\le \nw(X)\le\w$ and hence is metrizable. By Corollary~\ref{10.10}, the compact metrizable space $K$ is zero-dimensional and hence is homeomorphic to a closed subset of the Cantor cube $2^\w$ by Proposition~\ref{p:emb}(1). On the other hand, by Proposition~\ref{p:hur}(2), the uncountable compact metrizable space $K$ contains a closed subset homeomorphic to $2^\w$. Now we can apply Theorem~\ref{bernstein} and conclude that the space $K$ is weakly homeomorphic to the Cantor cube $2^\w$.
\smallskip

1.3. The space $X$ is not $\sigma$-compact. In this case we shall prove that $X$ is weakly homeomorphic to the Baire space $\IN^\w$. By Propositions~\ref{p:emb}(3) and \ref{p:hur}(3), $X$ is homeomorphic to a closed subset of $\IN^\w$ and contains a closed subspace homeomorphic to $\IN^\w$. By Theorem~\ref{bernstein}, the spaces $X$ and $\IN^\w$ are weakly homeomorphic.
\smallskip

2. Next assume that the space $X$ is $\sigma$-compact. If $X$ is Polish, then $X$ is homeomorphic to $\w$ or $2^\w$ by the items 1.1 and 1.2. So, we assume that the space $X$ is not Polish. Three subcases are possible.
\smallskip

2.1. The space $X$ is countable. Let $U$ be the union of all Polish open subspaces in $X$. It can be shown that the space $U$ is Polish and hence its complement $Z=X\setminus U$ is non-empty and nowhere Polish. In particular, $Z$ has no isolated points. By Proposition~\ref{p:char}(1), the countable space $Z$ is homeomorphic to $\IQ$, which implies that $X$ contains a closed subspace homeomorphic to $\IQ$. On the other hand, Proposition~\ref{p:emb}(1) guarantees that $X$ is homeomorphic to a closed subspace of $\IQ$. By Theorem~\ref{bernstein}, the spaces $X$ and $\IQ$ are weakly homeomorphic.
\smallskip

2.2. The space $X$ is contains a closed subspace homeomorphic to $\IQ\times 2^\w$. By Proposition~\ref{p:emb}(4), the space $X$ is homeomorphic to a closed subspace of $\IQ\times 2^\w$. By Theorem~\ref{bernstein}, the spaces $X$ and $\IQ\times 2^\w$ are weakly homeomorphic.
\smallskip

2.3. The space $X$ contains no closed copies of the space $\IQ\times 2^\w$. In this case we shall show that $X$ is weakly homeomorphic to $\IQ\oplus 2^\w$. By the decomposition Lemma~\ref{l:decomp}, the space $X$ is weakly homeomorphic to the topological sum $C\oplus P\oplus Z$ where $C$ is a metrizable locally countable space, $P$ is a metrizable locally Polish space and $Z$ is a closed subspace of $X$ such that every non-empty subset of $Z$ is not countable and not Polish. If $Z\ne\emptyset$, then by Proposition~\ref{p:char}(4), $Z$ is homeomorphic to $\IQ\times 2^\w$ and hence $X$ contains a closed copy of $\IQ\times 2^\w$, which contradicts our assumption.  Therefore, $Z=\emptyset$ and $X$ is weakly homeomorphic to the direct sum $C\oplus P$ of a locally countable and locally Polish spaces.
By Corollary~\ref{5.5}, $\hl(C\oplus P)\le \w\cdot\hl(X)\le\w$, which implies that the spaces $C,P$ are Lindel\"of and hence $C$ is countable and $P$ is Polish.
By Corollary~\ref{10.10}, the space $C\oplus P$ is $\sigma$-compact, which implies that the Polish space $P$ is $\sigma$-compact.

By Corollary~\ref{polish}, the space $C\oplus P$ is not Polish and hence the countable space $C$ is not Polish. By the case 2.1, the countable non-Polish space $C$ is weakly homeomorphic to $\IQ$. Since the space $X$ is uncountable and weakly homeomorphic to $C\oplus P$,  the Polish $\sigma$-compact space $P$ is uncountable and hence is weakly homeomorphic to $2^\w$ by the case 1.2. Consequently, the topological sum $C\oplus P$ is weakly homeomorphic to $\IQ\oplus 2^\w$. By Lemma~\ref{l:comp}, the space $X$ is weakly homeomorphic to $\IQ\oplus 2^\w$.
\smallskip

3. Finally assume that the space $X$ is $\sigma$-Polish. If $X$ is Polish or $\sigma$-compact, then by the preceding cases, $X$ is weakly homeomorphic to one of 6 spaces: $\w$, $2^\w$, $\IN^\w$, $\IQ$, $\IQ\times 2^\w$, $\IQ\oplus 2^\w$. So, we assume that $X$ is neither Polish nor $\sigma$-compact.
Two subcases are possible.
\smallskip

3.1. The space $X$ contains a closed subspace homeomorphic to $\IQ\times \IN^\w$.
By Proposition~\ref{p:emb}(5), $X$ is homeomorphic to a closed subspace of $\IQ\times\IN^\w$ and by Theorem~\ref{bernstein} the spaces $X$ and $\IQ\times\IN^\w$ are weakly homeomorphic.
\smallskip

3.2. The space $X$ contains no closed subspaces homeomorphic to $\IQ\times \IN^\w$.
By the decomposition Lemma~\ref{l:decomp}, the space $X$ is weakly homeomorphic to the topological sum $S\oplus P\oplus Z$ where $S$ is a metrizable locally $\sigma$-compact space, $P$ is a metrizable locally Polish space and $Z$ is a closed subspace of $X$ such that every non-empty subset of $Z$ is not $\sigma$-compact and not Polish. If $Z\ne\emptyset$, then by Proposition~\ref{p:char}(5), $Z$ is homeomorphic to $\IQ\times\IN^\w$ which contradicts our assumption. So, $Z=\emptyset$ and hence $X$ is weakly homeomorphic to $S\oplus P$. By Corollary~\ref{10.10}(2), the space $S\oplus P$ is Lindel\"of. Consequently, the locally $\sigma$-compact Lindel\"of space $S$ is $\sigma$-compact and the locally Polish Lindel\"of space $P$ is Polish. By Corollary~\ref{10.10}(4), the space $S\oplus P$ is zero-dimensional. By Corollaries~\ref{10.10}(6) and \ref{polish}, the space $S\oplus P$ is not $\sigma$-compact and not Polish. Consequently, the $\sigma$-compact space $S$ is not Polish and the Polish space $P$ is not $\sigma$-compact. By the case 1.3, the Polish zero-dimensional space $P$ is weakly homeomorphic to $\IN^\w$ and by case 2, the $\sigma$-compact non-Polish space $S$ is weakly homeomorphic to $\IQ$, $\IQ\oplus 2^\w$ or $\IQ\times 2^\w$. Then the space $S\oplus P$ is weakly homeomorphic to $\IQ\oplus\IN^\w$, $\IQ\oplus 2^\w\oplus\IN^\w$ or $(\IQ\times 2^\w)\oplus\IN^\w$. By Theorem~\ref{bernstein}, the space  $\IQ\oplus 2^\w\oplus\IN^\w$ is weakly homeomorphic to $\IQ\oplus\IN^\w$. Consequently, $X$ is weakly homeomorphic to $\IQ\oplus\IN^\w$ or $(\IQ\times 2^\w)\oplus\IN^\w$.
\end{proof}

\end{document}